\documentclass[11pt]{article}
\usepackage{amsmath,amsfonts,amsthm}
\usepackage[colorlinks]{hyperref}
\usepackage{breakurl}

\title{The Brauer group is not a derived invariant}
\author{Nicolas Addington}
\date{}

\renewcommand \O {\mathcal O}
\newcommand \Q {\mathbb Q}
\newcommand \Z {\mathbb Z}

\DeclareMathOperator \Br {Br}
\DeclareMathOperator \Hom {Hom}
\newcommand \Caldararu {C\u{a}l\-d\u{a}\-ra\-ru}
\DeclareMathOperator \tors {tors}
\renewcommand \top {\textrm{top}}
\DeclareMathOperator \Sq {Sq}

\newtheorem*{prop*}{Proposition}
\theoremstyle{definition}
\newtheorem*{def*}{Definition}

\begin{document}
\maketitle

\begin{abstract}
In this short note we observe that the recent examples of derived-equivalent Calabi--Yau 3-folds with different fundamental groups also have different Brauer groups, using a little topological K-theory.
\end{abstract}

Some years ago Gross and Popescu \cite{grpo} studied a simply-connected Calabi--Yau 3-fold $X$ fibered in non-principally polarized abelian surfaces.  They expected that its derived category would be equivalent to that of the dual abelian fibration $Y$, which is again a Calabi--Yau 3-fold but with $\pi_1(Y) = (\Z_8)^2$, the largest known fundamental group of any Calabi--Yau 3-fold.  This derived equivalence was later proved by Bak \cite{bak} and Schnell \cite{schnell}.  Ignoring the singular fibers, it is just a family version of Mukai's classic derived equivalence between an abelian variety and its dual \cite{mukai}, but of course the singular fibers require much more work.  As Schnell pointed out, it is a bit surprising to have derived-equivalent Calabi--Yau 3-folds with different fundamental groups, since for example the Hodge numbers of a 3-fold are derived invariants \cite[Cor.~C]{ps}.

Gross and Pavanelli \cite{grpa} showed that $\Br(X) = (\Z_8)^2$, the largest known Brauer group of any Calabi--Yau 3-fold.  In this note we will show that a certain extension of $\Br(X)$ by $H_1(X,\Z)$ is a derived invariant of Calabi--Yau 3-folds; thus in this example we must have $\Br(Y) = 0$, and in particular the Brauer group alone is not a derived invariant.  This too is a bit surprising, since the Brauer group \emph{is} a derived invariant of K3 surfaces: if $X$ is a K3 surface then $\Br(X) \cong \Hom(T(X), \Q/\Z)$ \cite[Lem.~5.4.1]{andrei}, where $T(X) = NS(X)^\perp \subset H^2(X,\Z)$ is the transcendental lattice, which is a derived invariant by work of Orlov \cite{orlov}.

After the first version this note circulated, Hosono and Takagi \cite{ht_der_eq} found a second example of derived-equivalent Calabi--Yau 3-folds with different fundamental groups.  Their $X$ and $Y$ are constructed from spaces of $5 \times 5$ symmetric matrices in what is likely an instance of homological projective duality \cite{ht_hpd}.  They satisfy $\pi_1(X) = \Z_2$ and $\pi_1(Y) = 0$, and while $\Br(X)$ and $\Br(Y)$ are not known, our result will give an exact sequence
\[ 0 \to \Z_2 \to \Br(Y) \to \Br(X) \to 0. \]
An explicit order-2 element of $\Br(Y)$ arises naturally in Hosono and Takagi's construction \cite[Prop.~3.2.1]{ht_der_eq}; presumably it is the image of $1 \in \Z_2$ above.

It is worth mentioning that both $\pi_1$ and $\Br$ are birational invariants, so while birational Calabi--Yau 3-folds are derived equivalent \cite{bridgeland}, the converse is not true.  In addition to the two examples just mentioned, there is the Pfaffian--Grassmannian derived equivalence of Borisov and \Caldararu\ \cite{bc}.  In that example $X$ is a complete intersection in a Grassmannian, so $H_1(X,\Z) = \Br(X) = 0$, so our result shows that $H_1(Y,\Z) = \Br(Y) = 0$ as well; to show that $X$ and $Y$ are not birational, Borisov and \Caldararu\ use a more sophisticated minimal model program argument.

Before proving our result we fix terminology.
\begin{def*}
A \emph{Calabi--Yau 3-fold} is a smooth complex projective 3-fold $X$ with $\omega_X \cong \O_X$ and $b_1(X) = 0$.  In particular $H_1(X,\Z)$ may be torsion.
\end{def*}
\noindent This is in contrast to the case of surfaces, where $\omega_X \cong \O_X$ and $b_1(X) = 0$ force $\pi_1(X) = 0$ \cite[Thm.~13]{kodaira}.  There are several reasons not to require $\pi_1(X) = 0$ for Calabi--Yau 3-folds.  As we have just seen, a simply-connected Calabi--Yau 3-fold may be derived equivalent to a non-simply-connected one; it may also be mirror to a non-simply-connected one.  Perhaps the best reason is that families of simply-connected and non-simply-connected Calabi--Yau 3-folds can be connected by ``extremal transitions,'' that is, by performing a birational contraction and then smoothing; most known families of Calabi--Yau 3-folds can be connected by extremal transitions \cite{gh,ks}.

\begin{def*}
The \emph{Brauer group} of a smooth complex projective variety $X$ is
\[ \Br(X) = \tors(H^2_\textrm{an}(X, \O_X^*)), \]
where $\tors$ denotes the torsion subgroup.
\end{def*}
\noindent This used to be called the \emph{cohomological Brauer group} until it was shown to coincide with the honest Brauer group \cite{djg}.
From the exact sequence
\[ H^2(X, \O_X) \to H^2(X, \O_X^*) \to H^3(X,\Z) \to H^3(X, \O_X) \]
we see that if $X$ is a Calabi--Yau 3-fold then
\[ \Br(X) = \tors(H^3(X,\Z)). \]
That is, the Brauer group of a Calabi--Yau 3-fold is entirely topological, in contrast to that of a K3 surface which is entirely analytic.

\begin{prop*}
Let $X$ and $Y$ be Calabi--Yau 3-folds with $D^b(X) \cong D^b(Y)$.  Then
\[ \left| H_1(X,\Z) \right| \cdot \left| \Br(X) \right| = \left| H_1(Y,\Z) \right| \cdot \left| \Br(Y) \right|. \]
More precisely, there is an exact sequence
\begin{equation} \label{K^1}
0 \to H_1(X,\Z) \to \tors(K^1_\top(X)) \to \Br(X) \to 0,
\end{equation}
where $K^*_\top$ is topological K-theory, and a similar sequence with $Y$; and an equivalence $D^b(X) \cong D^b(Y)$ induces an isomorphism $K^*_\top(X) \cong K^*_\top(Y)$.
\end{prop*}
\begin{proof}
Brunner and Distler \cite[\S2.5]{bd} analyzed the boundary maps in the Atiyah--Hirzebruch spectral sequence and saw that for a Calabi--Yau 3-fold $X$, or indeed any closed oriented 6-manifold with $b_1(X) = 0$, it degenerates at the $E_2$ page.  Thus there is an exact sequence
\[ 0 \to H^5(X,\Z) \to K^1_\top(X) \to H^3(X,\Z) \to 0. \]
Since $H^5(X,\Z) = H_1(X,\Z)$ is torsion, this yields the exact sequence \eqref{K^1}.\footnote{The published version of this paper cited \cite[\S4]{dm} to say that the sequence \eqref{K^1} splits.  But David Treumann has drawn my attention to several problems with this reference, especially the crucial \cite[Lem.~4.2]{dm}.  For one, the map $(c_1,c_2,c_3)$ is not a group homomorphism.  For another, on the quintic 3-fold, the class $[\O(1)] - [\O] \in \tilde K^0_\top$ satisfies $c_1^3 = 5$, $c_2 = 0$, $c_3 = 0$, hence is a counterexample to $\Sq^2 c_2 = [c_3 + c_1 c_2 + c_1^3]_2$.}
%
%
%
%
%

The fact that $K^*_\top$ is a derived invariant is discussed in \cite[\S2.1]{at}.  In a bit more detail, if $\Phi\colon D^b(X) \to D^b(Y)$ and $\Psi\colon D^b(Y) \to D^b(X)$ are inverse equivalences, then by \cite[Thm.~2.2]{orlov} there are objects $E, F \in D^b(X \times Y)$ such that
\begin{align*}
\Phi(-) &= \pi_{Y*}(E \otimes \pi_X^*(-)) &
\Psi(-) &= \pi_{X*}(F \otimes \pi_Y^*(-)),
\end{align*}
and arguing as in \cite[Lem.~5.32]{huybrechts} we find that the same formulas define inverse isomorphisms $K^*_\top(X) \to K^*_\top(Y)$ and $K^*_\top(Y) \to K^*_\top(X)$: use the fact that the pushforward on $K^*_\top$ satisfies a projection formula and is compatible with the pushforward on $D^b$.
\end{proof}

We conclude with a remark on $H_1$ and $\Br$ in mirror symmetry.  Batyrev and Kreuzer \cite{bk} predicted that mirror symmetry exchanges $H_1$ and $\Br$, having calculated both groups for all Calabi--Yau hypersurfaces in 4-dimensional toric varieties.  In all their examples the groups are quite small: either $H_1 = 0$ and $\Br = \Z_2$, $\Z_3$, or $\Z_5$, or vice versa.  This prediction does not seem to be right in general.  On the one hand it is contradicted by a prediction of Gross and Pavanelli \cite[Rem.~1.5]{grpa}, based on calculations in Pavanelli's thesis \cite{pavanelli}, that if $X$ is the abelian fibration above, with $H_1(X) = 0$ and $\Br(X) = (\Z_8)^2$, then its mirror $\check X$ has $\pi_1(\check X) = \Br(\check X) = \Z_8$.  Even more seriously, Hosono and Takagi's $X$ and $Y$ have the same mirror according to \cite{ht_mirror}, but different $H_1$ and $\Br$ as we have discussed.  Mirror symmetry is expected to exchange $K^0_\top$ and $K^1_\top$, however, and for a complex manifold (or more generally a $\operatorname{Spin}^c$ manifold) these have isomorphic torsion subgroups.
\vspace \baselineskip

I thank Paul Aspinwall and Andrei \Caldararu\ for helpful conversations, Shinobu Hosono and Hiromichi Takagi for encouraging me to publish this note, and David Treumann for an extensive discussion about the splitting of the exact sequence \eqref{K^1}.

\bibliographystyle{plain}
\bibliography{br_not_inv}

\noindent Nicolas Addington \\
Department of Mathematics \\
University of Oregon \\
Eugene, OR 97403-1222 \\[\baselineskip]
adding@uoregon.edu

\end{document}